\documentclass[twoside,10pt]{article}
\usepackage[]{amsmath,amsthm,amssymb}
\usepackage[latin1]{inputenc}

\begin{document}

\begin{center}{{\Large\bf  On the curious series related to the elliptic integrals }}\end{center}
\vspace{0,5cm}
\begin{center}{{\bf Semyon  Yakubovich}}\end{center}

\markboth{\rm \centerline{ Semyon  Yakubovich}}{}
\markright{\rm \centerline{Series related to the elliptic integrals}}

\begin{abstract} {\noindent By using the theory of the elliptic integrals a new method of summation is proposed for a certain class of series and their derivatives involving hyperbolic functions.   It is based on the termwise differentiation of the series with respect to the elliptic modulus and integral representations of several of the series in terms of the inverse Mellin transforms related to the Riemann zeta function.  The relation with the corresponding case of the Voronoi summation formula is exhibited. The involved series are expressed in closed form in terms of complete elliptic integrals of the first and second kind,  and some special  cases are calculated in terms of particular values of the Euler gamma function.}

\end{abstract}
\vspace{4mm}
{\bf Keywords}: {\it Series with hyperbolic functions, Elliptic integrals,  Mellin transform,
Riemann zeta function, Euler gamma function, arithmetic functions, summation formulae}

\vspace{2mm}

 {\bf AMS subject classification}:  40A99, 44A15, 33C75,  11K65, 11M06  

\vspace{4mm}

\section {Introduction and auxiliary results}

The main goal of this paper is to find closed-form relations  for the following series, involving the hyperbolic functions

$$\sum_{n=1}^\infty  (\pm 1)^n  n^\alpha {\cosh^\beta (\pi a n x) \over \sinh^\gamma (\pi an x)},\quad  \sum_{n=1}^\infty  (\pm 1)^n  n^\alpha {\sinh^\beta (\pi a n x) \over \cosh^\gamma (\pi an x)},\eqno(1.1)$$

$$\sum_{n=1}^\infty  n^\alpha \left[ \tanh(\pi a n x) -1\right],\quad \sum_{n=1}^\infty  n^\alpha  \left[ \coth(\pi a n x) -1\right],\eqno(1.2)$$
where $\alpha, \beta,\gamma  \in \{0,1, 2, 3\},\  \gamma > \beta,  \  a \in \{ 1/2, 1, 2 \}$,  $x >0$, being defined by the formula
$$x\equiv x(k)= { K (k^\prime)\over K(k)},\  k \in (0, 1),\ k^\prime = \sqrt {1-k^2},\eqno(1.3)$$
and $K(k)$ is the complete elliptic integral of the first kind [1], [7], Vol. II, [11]

$$K(k)= \int_0^1 {dt\over \sqrt{ (1-t^2) (1- k^2 t^2)} }.\eqno(1.4)$$  
The parameter $k$ is called the elliptic modulus and $k^\prime$ is the complementary modulus.   K(k) satisfies the Legendre relation

$$E(k)K (k^\prime) + E (k^\prime)K(k)-  K (k^\prime) K(k) = {\pi\over 2},\eqno(1.5)$$
where $E(k)$ is the complete elliptic integral of the second kind

$$E(k)= \int_0^1 \sqrt{{ 1- k^2 t^2\over 1-t^2} } dt.\eqno(1.6)$$  
Its derivative can be calculated by formula 

$${dE\over dk} = {E(k)- K(k)\over k}.\eqno(1.7)$$

It is known [11],  that $K(k),  K(k^\prime)$ satisfy  the differential equation

$${d\over dk } \left( k (k^\prime)^2 {du\over dk} \right)= k u \eqno(1.8)$$
and $E(k),\   E(k^\prime)- K(k^\prime)$ are, in turn,  solutions of  the differential equation

$$(k^\prime)^2 {d\over dk } \left( k  {du\over dk}\right) + ku =0.\eqno(1.9)$$  
The derivative of $K(k)$ can be calculated by the formula

$${dK\over dk} =   {E(k)-   (k^\prime)^2 K(k)\over k (k^\prime)^2} .\eqno(1.10)$$

In order to express series (1.1), (1.2) in closed form we will proposed a method of termwise differentiation with respect to the elliptic modulus and in some cases we will represent these series, using the inverse Mellin transform  related to the Riemann zeta function.  We note that this investigation is a continuation of earlier methods proposed by Ling and Zucker (see [6] and [15], respectively),  involving Weierstrassian and Jacobian elliptic functions (cf. [1]) and by Glasser et al. [4], basing on the Poisson summation formula.  It is also worth to mention that recently some infinite series of the Eisenstein type involving the hyperbolic functions were investigated in [5].   

Let $k_r$ be an elliptic modulus such that $x(k_r) = \sqrt r$ (see (1.3)).   In the sequel we will use such values for small $r$ and the corresponding elliptic integral singular values $K(k_r)$ (see [2], [3] ), namely

$$k_1= {1\over \sqrt 2}, \  k_2= \sqrt 2- 1,\  k_3= {1\over 4} \sqrt 2 (\sqrt 3-1),\ k_4= 3- 2\sqrt 2,\eqno(1.11)$$

$$K(k_1) = {\Gamma^2(1/4)\over 4\sqrt \pi},\quad     K(k_2) = {(\sqrt 2+1)^{1/2} \Gamma(1/8) \Gamma(3/8)\over 2^{13/4} \sqrt \pi},\eqno(1.12)$$

$$K(k_3) = {3^{1/4} \Gamma^3(1/3)\over 2^{7/3} \pi},\quad     K(k_4) = {(\sqrt 2+1) \Gamma^2 (1/4) \over 2^{7/2} \sqrt \pi}, \eqno(1.13)$$
where $\Gamma (z)$ is Euler's gamma function. According to [2] the so-called elliptic alpha function for the integral singular values 
$$\alpha(r)=  {E (k^\prime_r)\over K(k_r)} - {\pi\over 4 [K(k_r)]^2}=  {\pi\over 4 [K(k_r)]^2} + \sqrt r \left[ 1- 
{E (k_r)\over K(k_r)} \right] \eqno(1.14)$$
is calculated, in particular,  for small values and we have

$$\alpha(1)= {1\over 2},\   \alpha(2)=  \sqrt 2-1, \ \alpha(3)=  {1\over 2} (\sqrt 3-1),\  \alpha(4)=  2 (\sqrt 2-1)^2.\eqno(1.15)$$

Meanwhile,  appealing to relations (2.4.3.1), (2.4.3.3), (2.4.9.2)  in [9], Vol. I and the inverse Mellin transform [10], we derive the following integral representations, which will be useful in the sequel

$${1\over \sinh( cx)} = {1\over \pi i} \int_{\mu-i\infty}^{\mu+i\infty}  (2^s-1) \Gamma(s) \zeta(s) (2cx)^{-s} ds,\ c >0, \mu > 1,\eqno(1.16)$$

$${1\over \sinh^2 ( cx)} = {2\over \pi i} \int_{\mu-i\infty}^{\mu +i\infty}  \Gamma(s) \zeta(s-1) (2cx)^{-s} ds, \ c >0, \mu > 2, \eqno(1.17)$$

$${1\over \cosh^2( cx)} = {2\over \pi i} \int_{\mu-i\infty}^{\mu+i\infty}  (1- 2^{2-s}) \Gamma(s) \zeta(s-1 ) (2cx)^{-s} ds,\ c >0, \mu > 0, \eqno(1.18)$$

$$ \tanh( cx)- 1 = {1\over \pi i} \int_{\mu-i\infty}^{\mu+i\infty}  (2^{1-s} -1) \Gamma(s) \zeta(s) (2cx)^{-s} ds,\ c >0, \mu > 0, \eqno(1.19)$$

$$ \coth( cx)- 1 = {1\over \pi i} \int_{\mu-i\infty}^{\mu+i\infty}  \Gamma(s) \zeta(s) (2cx)^{-s} ds,\ c >0, \mu > 1, \eqno(1.20)$$
where $\zeta(s)$ is the Riemann zeta function [7], Vol. I, which satisfies the  familiar functional equation
$$\zeta(s)= 2^s \pi^{s-1} \sin\left({\pi s\over 2}\right)\Gamma(1-s)\zeta(1-s).\eqno(1.21)$$

\section{Series with  the hyperbolic functions}

In order to evaluate series (1.1) we propose the method of termwise differentiation of the series with respect to the elliptic modulus.  To do this,  we will employ  formulae of Sections 5.3.4. and 5.3.6. in  [9], Vol. I.   Indeed, let us consider the following series (see relation (5.3.4.2) in [9], Vol. I)

$$\sum_{n=1}^\infty  {n \over \sinh(\pi n x)} = {1\over \pi^2} K(k) \left[ K(k)- E(k) \right],\eqno(2.1)$$
where $x(k) $ is defined by (1.3).  It is easy to verify that the function $x: (0,1) \to \mathbb{R}_+$ is monotone decreasing and continuously differentiable.  This means that equality (2.1) is valid for any $x >0$.  Moreover,  the series (2.1) is differentiated  termwise with respect to $k \in (0, k_0],\ k_0 < 1$ via the absolute and uniform convergence of its derivative.   Hence we have

{\bf Theorem 1}. {\it     For all $x >0$ the following identities hold}

$$  2  \sum_{n=1}^\infty  {n^2 \cosh(\pi n x)  \   \over  \sinh^2 (\pi n x)}  =   \sum_{n=1}^\infty  {(2n-1) \cosh( \pi (2 n-1)  x/2)  \over \sinh^3( \pi (2 n-1)  x/2)}$$

$$=  {4\over \pi^4} K^2(k)   \left[ k^2 K^2(k) -   \left[ K(k)- E(k)\right]^2\right]. \eqno(2.2)$$

\begin{proof}   In fact, as we see above $x(k)$ is a bijective map from $(0,1)$ to $ \mathbb{R}_+$. Moreover,  termwise differentiation with respect to $k$ in (2.1)  gives 

$$- \pi x^\prime(k) \sum_{n=1}^\infty  {n^2 \cosh(\pi n x)  \over \sinh^2 (\pi n x)} = {1\over \pi^2} {d\over dk} \left[K(k) \left[ K(k)- E(k) \right] \right].\eqno(2.3)$$
Meanwhile, since via relation (5.3.4.6) in [9], Vol. I

$$\sum_{n=1}^\infty  {1 \over \sinh^2 (\pi (2n-1) x/2)} = {2\over \pi^2} K(k) \left[ K(k)- E(k) \right],\eqno(2.4)$$
then after differentiation we easily establish the first equality in (2.2).  Further,  employing twice (1.10), we find

$$x^\prime (k) =  -  \frac{K(k^\prime)} {K^2(k)} {dK(k)\over dk} -  \frac{k } {(1-k^2) K(k)} \left[ \frac{E (k^\prime)} {k^2} -
K(k^\prime)\right]$$

$$=     \frac{1} {k(1-k^2) K(k)} \left[ K(k^\prime) \left[ 1 -  \frac{E(k)} {K(k)} \right] -   E (k^\prime)\right] $$ 
and the Legendre identity (1.5) leads  us  to the final result

$$x^\prime (k) =  - \frac{\pi}{2 k(1-k^2) K^2(k)}.\eqno(2.5)$$
On the other hand,  with the aid of (1.7), (1.9)  and (1.10)

$${d\over dk} \left[K(k) \left[ K(k)- E(k) \right] \right] =  {dK(k) \over dk}  \left[ K(k)- E(k) \right] 
-    K(k) {d\over dk } \left( k {dE(k)\over dk} \right) $$

$$=    \left[ K(k)- E(k) \right]  {E(k)-   (1- k^2)  K(k)\over k (1- k^2) }  + {k\  K(k) E(k) \over 1- k^2}  $$

$$=   {k\ K^2(k) \over 1- k^2}  -    {\left[ K(k)- E(k)\right]^2  \over k(1- k^2)}. $$ 
Hence, (2.3) and (2.5) yield the equality 

$$ \sum_{n=1}^\infty  {n^2 \cosh(\pi n x)  \over \sinh^2 (\pi n x)} = {2\over \pi^4} K^2(k) 
\left[ k^2 K^2(k) -   \left[ K(k)- E(k)\right]^2\right],$$
which implies the latter equality in (2.2) and completes the proof of Theorem 1.
\end{proof}

Appealing to (1.11), (1.12), (1.13),  (1.15),   we arrive at an immediate corollary.

{\bf Corollary 1}.    {\it  The following formulae take place}

$$ 2  \sum_{n=1}^\infty  {n^2 \cosh(\pi n)  \   \over  \sinh^2 (\pi n)}  =   \sum_{n=1}^\infty  {(2n-1) \cosh( \pi (n- 1/2) ) \over \sinh^3( \pi ( n-1/2))}$$

$$=  {1\over 4 \pi^2}  \left[ {\Gamma^8(1/4) \over 64  \pi^4}  +  {\Gamma^4 (1/4) \over 4  \pi^2} - 1\right],\eqno(2.6)$$

$$ 2  \sum_{n=1}^\infty  {n^2 \cosh(\pi n \sqrt 2)  \   \over  \sinh^2 (\pi n \sqrt 2)}  =   \sum_{n=1}^\infty  {(2n-1) \cosh( \pi (2n- 1)/\sqrt 2)  \over \sinh^3( \pi (2 n-1)/\sqrt 2)}$$

$$=  {1\over 8 \pi^2}  \left[ {\Gamma^4(1/8) \Gamma^4(3/8)   \over 512  \pi^4}  +  {\Gamma^2 (1/8) \Gamma^2(3/8) \over 8\sqrt 2   \pi^2} - 1\right],\eqno(2.7)$$

$$ 2  \sum_{n=1}^\infty  {n^2 \cosh(\pi n \sqrt 3)  \   \over  \sinh^2 (\pi n \sqrt 3)}  =   \sum_{n=1}^\infty  {(2n-1) \cosh( \pi (n- 1/2)\sqrt 3 )  \over \sinh^3( \pi ( n-1/2)\sqrt 3)}$$

$$=  {1\over 12 \pi^2}  \left[ {3\left( 2^{- 1/3}\right)   \Gamma^{12} (1/3) (2-\sqrt 3)  \over 128 \pi^6  }  +  {\sqrt 3 \left(2^{-2/3}\right)   \Gamma^6 (1/3) (\sqrt 3- 1)  \over  4 \pi^3} - 1\right],\eqno(2.8)$$

$$ 2  \sum_{n=1}^\infty  {n^2 \cosh(2\pi n )  \   \over  \sinh^2 (2\pi n )}  =   \sum_{n=1}^\infty  {(2n-1) \cosh( \pi (2n- 1) )  \over \sinh^3( \pi (2 n-1))}$$

$$=  {1\over 16 \pi^2}  \left[ { \Gamma^{4} (1/4)  \over 8 \pi^2  }   - 1\right].\eqno(2.9)$$

Next, we will combine relations (5.3.4.1) and (5.3.6.2) in [9], Vol. I, namely,

$$\sum_{n=1}^\infty  {(-1)^{n-1}  \over \sinh(\pi (2n-1) x/2)} = {k\over \pi} K(k),\eqno(2.10)$$

$$\sum_{n=1}^\infty  {1 \over \cosh(\pi nx )} = {1\over \pi} K(k)-  {1\over 2},\eqno(2.11)$$
to derive other new formulae.   In fact, termwise differentiation with respect to $k$ and the use of (1.10), (2.5) give the results

$$\sum_{n=1}^\infty  {(-1)^{n-1} (2n-1) \cosh (\pi (2n-1) x/2) \over \sinh^2(\pi (2n-1) x/2)} = {4 k \over \pi^3}\   K^2(k) E(k),\ x >0,\eqno(2.12)$$

$$ \sum_{n=1}^\infty  {n \sinh(\pi nx)  \over \cosh^2 (\pi nx )} = {2\over \pi^3}   K^2(k) \left[ E(k)-   (1-k^2) K(k) \right],\ x >0.\eqno(2.13)$$
Obtaining particular cases in the same fashion as in Corollary 1, we establish

{\bf Corollary 2}.     {\it  The following formulae take place}

$$\sum_{n=1}^\infty  {(-1)^{n-1} (2n-1) \cosh (\pi (n-1/2) )\over \sinh^2(\pi (n-1/2))} = {\Gamma^2(1/4) \over 4 \pi^2 \sqrt{2\pi}}\  \left[ {\Gamma^4(1/4) \over 8 \pi^2} + 1\right],\eqno(2.14)$$

$$ \sum_{n=1}^\infty  {n \sinh(\pi n)  \over \cosh^2 (\pi n )} =  { \Gamma^2(1/4) \over 8 \pi^2 \sqrt{\pi}} ,\eqno(2.15)$$

$$\sum_{n=1}^\infty  {(-1)^{n-1} (2n-1) \cosh (\pi (2n-1) /\sqrt 2) \over \sinh^2(\pi (2n-1)/\sqrt 2)} =
\frac{ 2^{1/4} \Gamma(1/8) \Gamma(3/8) } {16 \pi^2 \sqrt \pi} $$

$$\times \left[  \frac{ (\sqrt 2+1)^{1/2} \Gamma^2(1/8) \Gamma^2(3/8) } {16 \pi^2 \sqrt 2} + (\sqrt 2-1)^{1/2} \right],\eqno(2.16)$$

$$ \sum_{n=1}^\infty  {n \sinh(\pi n\sqrt 2)  \over \cosh^2 (\pi n\sqrt 2 )} =
\frac{  \Gamma(1/8) \Gamma(3/8) } {16 \pi^2 \sqrt {2\pi}} $$

$$\times \left[  \left( 1+ {1\over \sqrt 2} \right)^{1/2} - \   \frac{ (2- \sqrt 2)^{1/2} \Gamma^2(1/8) \Gamma^2(3/8) } {32 \pi^2 }  \right],\eqno(2.17)$$

$$\sum_{n=1}^\infty  {(-1)^{n-1} (2n-1) \cosh (\pi (n-1/2) \sqrt 3) \over \sinh^2(\pi (n-1/2) \sqrt 3)} $$

$$= \frac{  3^{1/4}\  \Gamma^3(1/3)  } {8\sqrt 2\  \pi^3}  \left[   \frac{ \Gamma^6(1/3) } {8 \pi^3 }  +  {1\over  2^{1/3} } \left( 1- {1\over \sqrt 3} \right)  \right],\eqno(2.18)$$

$$ \sum_{n=1}^\infty  {n \sinh(\pi n\sqrt 3)  \over \cosh^2 (\pi n\sqrt 3 )} =   \frac{  3^{- 1/4}\  \Gamma^3(1/3)  } {8 \  \pi^3}  \left[ 2^{-1/3} -    \frac{ \sqrt 3\  \Gamma^6(1/3) } {32   \pi^3 } \right],\eqno(2.19)$$ 

$$\sum_{n=1}^\infty  {(-1)^{n-1} (2n-1) \cosh (\pi (2n-1) ) \over \sinh^2(\pi (2n-1) )} = 
 {\Gamma^2(1/4)\over 16 \pi^2 \sqrt {2\pi} }  \left[ { \Gamma^{4} (1/4)  \over 8 \pi^2  }   + \sqrt 2- 1\right],\eqno(2.20)$$

$$ \sum_{n=1}^\infty  {n \sinh(2 \pi n)  \over \cosh^2 (2\pi n)} =  {\Gamma^2(1/4)\over 32  \pi^2 \sqrt {2\pi} }  \left[\sqrt 2 + 1-   { \Gamma^{4} (1/4)  \over 8   \pi^2  } \right].\eqno(2.21)$$

Taking relations (5.3.4.3), (5.3.4.4), (5.3.4.5),  (5.3.6.4), (5.3.6.5), (5.3.6.6) in [9], Vol. I

$$ \sum_{n=1}^\infty  {(-1)^{n-1} \  n \over \sinh(\pi n x)} = {K(k) \over \pi^2} \left[ E(k)-  (1-k^2) K(k)\right],\eqno(2.22)$$

$$ \sum_{n=1}^\infty  {1 \over (2n-1) \sinh(\pi (2n-1) x)} = - {1\over 8} \log (1-k^2),\eqno(2.23)$$

$$ \sum_{n=1}^\infty  {1 \over \sinh^2(\pi n x)} = {1\over 6} - {2 K(k) \over \pi^2} \left[ E(k) -  {2-k^2\over 3} K(k) \right],\eqno(2.24)$$

$$ \sum_{n=-\infty}^\infty  {1  \over \cosh (\pi (2n-1) x/2)} = {2 k\over \pi} K(k),\eqno(2.25)$$

$$ \sum_{n=1}^\infty  {1 \over \cosh^2(\pi n x)} =  {2\over \pi^2} K(k)E(k) - {1\over 2},\eqno(2.26)$$

$$ \sum_{n=1}^\infty  {1 \over \cosh^2(\pi (2n-1) x/2)} =  {2 K(k) \over \pi^2} \left[ E(k) -  (1-k^2) K(k)\right],\eqno(2.27)$$
respectively,   recalling (2.1), (2.4), (2.10)  and summing or subtracting one from another,  we deduce, for instance, the following equalities

$$ \sum_{n=-\infty}^\infty  {1  \over \cosh (\pi (2n-1) x/2)} +  \sum_{n=1}^\infty  {(-1)^n  \over \sinh (\pi (2n-1) x/2)}  = {k\over \pi} K(k),\eqno(2.28)$$

$$ \sum_{n=1}^\infty  { \cosh(2\pi nx)  \over \sinh^2 (2\pi n x)}  =   {2-k^2\over 6 \pi^2} K^2(k) - {1\over 12},\eqno(2.29)$$

$$ \sum_{n=1}^\infty  { 1 \over \sinh^2 (2\pi n x)}  =   {K(k) \over \pi^2} \left[ {2-k^2\over 6 } K(k) -  E(k)\right] +   {1\over 6},\eqno(2.30)$$

$$ \sum_{n=1}^\infty  { 1 \over \cosh^2 (\pi n x/2)}  =   {2 K(k) \over \pi^2} \left[ 2E(k) -    (1-k^2)  K(k) \right]  - {1\over 2},\eqno(2.31)$$

$$ \sum_{n=1}^\infty  { 1 \over \cosh^2 (\pi n x)}  -   \sum_{n=1}^\infty  {1 \over \cosh^2(\pi (2n-1) x/2)} =   {2(1-k^2) \over  \pi^2} K^2(k) - {1\over 2},\eqno(2.32)$$

$$ \sum_{n=1}^\infty  {(-1)^{n-1} \  n \over \sinh(\pi n x)} = \sum_{n=1}^\infty  {1 \over 2 \cosh^2(\pi (2n-1) x/2)},\eqno(2.33)$$

$$ \sum_{n=1}^\infty  {1 \over \sinh^2(\pi (2n-1) x)} =   {K(k)  \over \pi^2}  \left[ {2-k^2\over 2} K(k)  - E(k)\right] .\eqno(2.34)$$
Combining (2.34) with (2.4), we get 

$$\sum_{n=1}^\infty  {\cosh(\pi(2n-1)x) \over  \sinh^2(\pi (2n-1) x)} = {k^2\over 2\pi^2} K^2(k).\eqno(2.35)$$
We note that the same result can be obtained, differentiating (2.23) with respect to $k$ and invoking (2.29).  An immediate corollary of (2.29), (2.35) is the equality

$$\sum_{n=1}^\infty  { \cosh(\pi nx) \over \sinh^2 (\pi n x)}  =  {1+ k^2\over 3 \pi^2} K^2(k)   - {1\over 12},\eqno(2.36)$$

Now we are ready to apply  the method of termwise differentiation with respect to the elliptic modulus to the series  (2.22), (2.24),  (2.25),  (2.26),  (2.27), (2.28), (2.29), (2.36).   In fact,   employing  (2.5) and properties of the complete elliptic integrals listed in Section 1, in particular,  (1.7), (1.10), after elementary calculations we  establish

{\bf Theorem 2}. {\it For all $x >0$ the  following formulae hold valid}

$$ \sum_{n=1}^\infty  {(-1)^{n-1}\  n^2  \cosh(\pi n x) \over \sinh^2 (\pi n x)} = {2\over \pi^4} K^2(k) \left[ \left( E(k) - (1-k^2) K(k) \right)^2+  k^2(1-k^2) K^2(k) \right],\eqno(2.37)$$

$$ \sum_{n=1}^\infty  {n \cosh(\pi n x)  \over \sinh^3 (\pi n x)} = {2\over 3 \pi^4} K^2(k) \left[  E(k) \left( 2 (2-k^2)K(k) -3 E(k) \right) -  (1-k^2)K^2(k) \right],\eqno(2.38)$$

$$ \sum_{n=1}^\infty  {(2n-1) \sinh(\pi (2n-1) x/2)   \over \cosh^2 (\pi (2n-1) x/2)} = {4 k\over \pi^3 } K^2(k) E(k),\eqno(2.39)$$

$$ \sum_{n=1}^\infty  {n \sinh(\pi n x)  \over \cosh^3 (\pi n x)} =  {2\over \pi^4} K^2(k) \left[ E^2(k) - (1-k^2) K^2(k)\right],\eqno(2.40)$$

$$ \sum_{n=1}^\infty  {(2n-1)\sinh(\pi (2n-1) x/2)  \over \cosh^3 (\pi (2n-1) x/2)} =  {4\over \pi^4} K^2(k) \left[ \left( E (k) - (1-k^2)K(k) \right)^2 + k^2 (1-k^2) K^2(k)\right]$$

$$= {2\over \pi} \sum_{n=1}^\infty  {(-1)^n\  n^2  \cosh(\pi n x) \over \sinh^2 (\pi n x)},\eqno(2.41)$$

$$ \sum_{n=-\infty}^\infty  {(2n- 1) \sinh (\pi (2n-1) x/2)  \over \cosh^2 (\pi (2n-1) x/2)} +  \sum_{n=1}^\infty  {(-1)^n (2n-1) \cosh (\pi (2n-1) x/2)  \over \sinh^2 (\pi (2n-1) x/2)}$$

$$  = {4k \over \pi^3}   K^2(k) E(k),\eqno(2.42)$$

$$ \sum_{n=1}^\infty  { n (3+ \cosh(4\pi nx) ) \over \sinh^3 (2\pi n x)}  =  {2\over 3\pi^4} K^3(k) \left[ (2-k^2)  E(k) -  2 (1-k^2)K (k) \right],\eqno(2.43)$$

$$ \sum_{n=1}^\infty  { n (3+ \cosh(2\pi nx) ) \over \sinh^3 (\pi n x)}  =  {8 \over 3\pi^4} K^3(k) 
\left[ E(k) (1+k^2)- K(k) (1-k^2) \right],\eqno(2.44)$$

$$   \sum_{n=1}^\infty  {n  \cosh(2\pi n x)   \over \sinh^3 (2\pi n x)} = {K^2(k) \over 6 \pi^4}  \left[ (2-k^2) E(k)K(k)  +  (1-k^2)K^2 (k)- 3 E^2(k) \right].\eqno(2.45)$$

{\bf Remark 1}.  In the same manner similar equalities can be obtained, differentiating  (2.30), (2.31), (2.32), (2.33),  (2.34), (2.35) termwise  with respect to the elliptic modulus $k$.

Some particular values of the above series are listed in 

{\bf Corollary 3}.  {\it It has a set of identities}  

$$ \sum_{n=1}^\infty  {(-1)^{n-1} \  n^2  \cosh(\pi n) \over \sinh^2 (\pi n)} = {1\over 8\pi^2}\left[1   + {\Gamma^8(1/4) \over 64 \pi^4}\right],$$

$$ \sum_{n=1}^\infty  {(2n-1)\sinh(\pi (n-1/2))  \over \cosh^3 (\pi (n-1/2))} = {1\over 4\pi^2}\left[1   + {\Gamma^8(1/4) \over 64 \pi^4}\right],$$

$$ \sum_{n=1}^\infty  {n \cosh(\pi n)  \over \sinh^3 (\pi n)} = {1 \over  8 \pi^2}  \left[  { \Gamma^8(1/4) \over 192 \pi^4}  - 1 \right],$$

$$ \sum_{n=1}^\infty  {(2n-1) \sinh(\pi (n-1/2) )   \over \cosh^2 (\pi (n-1/2))} = { \Gamma^2 (1/4) \over  4 \pi^2 \sqrt {2\pi}  }  \left[ 1 + {\Gamma^4(1/4) \over 8 \pi^2 }\right],$$

$$ \sum_{n=1}^\infty  {n \sinh(\pi n)  \over \cosh^3 (\pi n)} =  {1\over 8\pi^2}  \left[ 1  +  {\Gamma^4(1/4)\over 4 \pi^2} -  {\Gamma^8(1/4)\over 64 \pi^4}\right],$$

$$\sum_{n=1}^\infty  { \cosh(\pi(2n-1)) \over \sinh^2 (\pi (2n-1))}  =   {\Gamma^4(1/4)\over 64 \pi^3},$$ 

$$\sum_{n=1}^\infty  { \cosh(\pi n) \over \sinh^2 (\pi n)}  =   {\Gamma^4(1/4)\over 32  \pi^3}-  {1\over 12},$$ 

$$ \sum_{n=-\infty}^\infty  {(2n- 1) \sinh (\pi (n- 1/2))  \over \cosh^2 (\pi (n-1/2))} +  \sum_{n=1}^\infty  {(-1)^n (2n-1) \cosh (\pi (n-1/2))  \over \sinh^2 (\pi (n-1/2)) }$$

$$  = { \Gamma^2(1/4) \over 4 \pi^3\  \sqrt {2\pi} }   \left[ \pi+  {\Gamma^4(1/4) \over 8 \pi } \right],$$

$$ \sum_{n=1}^\infty  { n (3+ \cosh(4\pi n) ) \over \sinh^3 (2\pi n )}  =   {\Gamma^4(1/4) \over 64  \pi^4 }
\left[ 1 -   {\Gamma^4(1/4) \over 24 \pi^2 } \right],$$

$$ \sum_{n=1}^\infty  { n (3+ \cosh(2\pi n) ) \over \sinh^3 (\pi n)}  =   {\Gamma^4(1/4) \over 48 \pi^4 }
\left[ 3 +  {\Gamma^4(1/4) \over 8 \pi^2 } \right],$$

$$   \sum_{n=1}^\infty  {n  \cosh(2\pi n)   \over \sinh^3 (2\pi n)} = {1 \over 32  \pi^2}  \left[   {\Gamma^8(1/4) \over 96  \pi^4 } -  {\Gamma^4(1/4)\over 8\pi^2} - 1  \right].$$

{\bf Remark 2}.  As a conclusion we stress that employing formulae from  Sections (5.3.4) and (5.3.6) in  [9], Vol. I, one can make different combinations of the known series  and differentiate them   termwise any number of times with respect to the elliptic modulus $k$ to obtain the values of new series, involving powers of hyperbolic functions.

In the meantime,  series (2.24), (2.26)  can be expressed using integral representations  (1.17), (1.18).  In fact, substituting the corresponding integrals inside the series, we  change  the order of integration and summation. This is  allowed by virtue of the absolute and uniform convergence with respect to $x \ge x_0 >0$, basing , in turn, on the estimate (it concerns series (2.24) and (2.26) can be treated analogously)

$$\sum_{n=1}^\infty {1\over \sinh^2(\pi nx)} = {2\over \pi} \sum_{n=1}^\infty \left|  \int_{\mu-i\infty}^{\mu+i\infty}  \Gamma(s) \zeta(s-1) (2\pi n x)^{-s} ds \right| $$

$$\le  2^{1-\mu}  \pi^{-\mu-1} x_0^{-\mu} \zeta(\mu-1) \zeta(\mu)  \int_{\mu-i\infty}^{\mu+i\infty}  \left| \Gamma(s) ds\right| < \infty, \quad \mu  > 2,$$
where the convergence of the latter integral can be verified using the Stirling asymptotic formula for the gamma-function when $|{\rm Im} s| \to \infty$ (see [7], Vol. I).   Hence, we find the representations 

$$ \sum_{n=1}^\infty {1\over \sinh^2(\pi nx)} =  {2\over \pi i} \int_{\mu-i\infty}^{\mu+i\infty}  \Gamma(s) \zeta(s) \zeta(s-1) (2\pi x)^{-s} ds,\eqno(2.46)$$

$$ \sum_{n=1}^\infty {1\over \cosh^2(\pi nx)} =  {2\over \pi i} \int_{\mu-i\infty}^{\mu+i\infty} \left(1- 2^{2-s}\right)  \Gamma(s) \zeta(s) \zeta(s-1) (2\pi x)^{-s} ds.\eqno(2.47)$$
Meanwhile, the product of the Riemann zeta-functions can be represented by the Ramanujan identity [13]

$$\zeta(s)\zeta (s-1) = \sum_{n=1}^\infty  {\sigma (n) \over n^s},\    {\rm Re} s > 2,\eqno(2.48)$$
where $\sigma(n)$ is the sum of the divisors of $n$ [1]. Therefore, substituting this expression in (2.46), (2.47) and changing the order of integration and summation due to the same motivation,  and employing the inverse Mellin transform of the gamma-function [10], we obtain the following equalities 

$$ \sum_{n=1}^\infty {1\over \sinh^2(\pi nx)} =  4 \sum_{n=1}^\infty \sigma(n) e^{-2\pi nx},\eqno(2.49)$$ 

$$ \sum_{n=1}^\infty {1\over \cosh^2(\pi nx)} =  4 \sum_{n=1}^\infty \sigma(n) \left[ e^{-2\pi nx} - 4\  e^{- 4\pi nx} \right].\eqno(2.50)$$ 
On the other hand, the Nasim summation formula [8] says

$$\sum_{n=1}^\infty \sigma(n) e^{-2\pi nx} + x^{-2} \sum_{n=1}^\infty \sigma(n) e^{-2\pi n/ x} = {1\over 24} \left( 1+ {1\over x^2}\right)- {1\over 4\pi x},\ x >0.$$
Hence,  appealing to (2.49), (2.50),  we establish the formulae  $(x >0)$ 

$$ \sum_{n=1}^\infty {1\over \sinh^2(\pi nx)} + {1\over x^{2}}   \sum_{n=1}^\infty {1\over \sinh^2(\pi n/ x)} = {1\over 6}\left(1 + {1\over x^2}\right) - {1\over \pi x},\eqno(2.51)$$

$$ \sum_{n=1}^\infty {1\over \cosh^2(\pi n x)} + {1 \over x^{2}}   \sum_{n=1}^\infty {1\over \cosh^2( \pi n/ x)} = 
{1\over 6} \left( 1+ {1\over x^2}\right)- {1\over \pi x}$$

$$- 4 \left(  \sum_{n=1}^\infty {1\over \sinh^2(2\pi nx)} + {1\over x^{2}}   \sum_{n=1}^\infty {1\over \sinh^2(2\pi n/ x)}\right).\eqno(2.52)$$

{\bf Theorem 3}. {\it  For all $x >0$  the following identities hold}

$$ \sum_{n=1}^\infty {1\over \sinh^2(\pi n/ x)} =  {1\over 6} -  {2\over \pi^2} K (k^\prime) \left[ E(k^\prime)
-  {k^2+1\over 3} K(k^\prime)\right],\eqno(2.53)$$

$$\sum_{n=1}^\infty {1\over \cosh^2(\pi n/ x)} =   {2\over \pi^2} K(k^\prime)E(k^\prime) - {1\over 2}.\eqno(2.54)$$

\begin{proof} Indeed, associating  any positive $x$ with some $k \in (0,1)$ by formula (1.3),   we appeal to (2.51) and (2.24) to deduce

$$\sum_{n=1}^\infty {1\over \sinh^2(\pi n/ x)} =  {1\over 6} - {x\over \pi} +  {2\over \pi^2} K^2 (k^\prime){E(k)\over K(k)}  -  {2(2-k^2)\over 3\pi^2} K^2(k^\prime).$$
But the Legendre identity (1.5) says

$${2\over \pi^2} K^2 (k^\prime){E(k)\over K(k)} =  {x\over \pi} +  {2\over \pi^2} K(k^\prime) \left[ K(k^\prime) - E(k^\prime) \right].$$ 
This drives us to (2.53).  In order to prove (2.54), we recall (2.4), (2.24) to calculate the series

$$\sum_{n=1}^\infty {1\over \sinh^2(\pi n x/2)}  =   \sum_{n=1}^\infty {1\over \sinh^2(\pi nx)} +  \sum_{n=1}^\infty {1\over \sinh^2(\pi (2n-1) x/2 )}$$

$$=   {1\over 6} - {4\over \pi^2} K(k)E(k) + {2(5-k^2)\over 3\pi^2} K^2(k).$$
Hence from (2.51), (1.5)  we get the value of the series

$$  \sum_{n=1}^\infty {1\over \sinh^2(2 \pi n/ x)} =  {1\over 6}  -  {K(k^\prime) \over \pi^2} \left[ E(k^\prime) -  {1+k^2\over 6} K(k^\prime)\right].$$
Consequently, employing identities (2.26), (2.30), (2.52), we establish (2.54), completing the proof of Theorem 3. 

\end{proof}

Summing (2.53), (2.54), we arrive at the immediate 

{\bf Corollary 4}. {\it It has the equality}

$$\sum_{n=1}^\infty {\cosh(2\pi n/ x)  \over \sinh^2(2\pi n/ x)} =   {1+k^2 \over 6 \pi^2}  K^2(k^\prime)  - {1\over 12}, \quad x >0.$$

In order to treat series (1.2),  we recall integral representations (1.19), (1.20)  and motivating the interchange of the order of integration and summation in the same manner as in (2.46), (2.47), we  use them to derive the identities 

$$\sum_{n=1}^\infty  n \left[ \tanh(\pi n x) -1\right] =   {1\over \pi i} \int_{\mu-i\infty}^{\mu+i\infty}  (2^{1-s} -1) \Gamma(s) \zeta(s)\zeta(s-1) (2\pi x)^{-s} ds $$

$$=  - {1\over 4}  \left[ \sum_{n=1}^\infty {1\over \cosh^2(\pi n x)}  +  \sum_{n=1}^\infty {1\over \sinh^2(\pi n x)} \right],\  x >0,$$

$$\sum_{n=1}^\infty  n \left[ \coth(\pi n x) -1\right] =   {1\over \pi i} \int_{\mu-i\infty}^{\mu+i\infty}  \Gamma(s) \zeta(s)\zeta(s-1) (2\pi x)^{-s} ds $$

$$=  {1\over 2}   \sum_{n=1}^\infty {1\over \sinh^2(\pi n x)},\  x >0.$$
Hence,   appealing to (2.4),  (2.24), (2.26), (2.29), (2.30), (2.38), (2.43), (2.45), (2.46)     and employing again the method of termwise differentiation  with respect to the elliptic modulus, we proved   the following results.

{\bf Theorem 4}.  {\it Let $ x>0$. Then }

$$\sum_{n=1}^\infty  n \left[1-  \tanh(\pi n x)\right] =  \sum_{n=1}^\infty {\cosh(2\pi n x) \over \sinh^2(2\pi n x)}= {1\over 12} \left[  {2(2-k^2)\over \pi^2} K^2(k) - 1 \right],\eqno(2.55)$$

$$\sum_{n=1}^\infty  n \left[ \coth(\pi n x) -1\right] =    {1\over 2}   \sum_{n=1}^\infty {1\over \sinh^2(\pi n x)}
= {1\over 12} - {K(k) \over \pi^2} \left[ E(k) -  {2-k^2\over 3} K(k)\right],\eqno(2.56)$$

$$\sum_{n=1}^\infty   {n^2 \over  \sinh^2(\pi n x)} =  {2 K^2(k) \over 3 \pi^4} \left[  E(k) \left( 2 (2-k^2)K(k) -3 E(k) \right) -  (1-k^2)K^2(k) \right],\eqno(2.57)$$

$$\sum_{n=1}^\infty  { n^2 \over  \cosh^2 (\pi n x)} = {2 K^3(k) \over 3\pi^4}  \left[ (2-k^2)  E(k) -  2 (1-k^2)K (k) \right].\eqno(2.58)$$

{\bf Theorem 5}.  {\it Let $ x>0$. Then }

$$\sum_{n=1}^\infty  n \left[ \coth(2\pi n x) -1\right] =   {1\over 12} - {K(k) \over 2\pi^2} \left[ E(k) -  {2-k^2\over 6} K(k)\right],\eqno(2.59)$$

$$\sum_{n=1}^\infty  n \left[1-  \tanh(\pi n x/2 )\right] =   {1\over 12} \left[  {4(1+k^2)\over \pi^2} K^2(k) - 1 \right],\eqno(2.60)$$

$$\sum_{n=1}^\infty   {n^2 \over  \sinh^2(2\pi n x)} =  {K^2(k) \over 6 \pi^4}  \left[ (2-k^2) E(k)K(k)  +  (1-k^2)K^2 (k)- 3 E^2(k) \right],\eqno(2.61)$$

{\bf Corollary 5}.    {\it The following relations hold}

$$\sum_{n=1}^\infty  n \left[1-  \tanh(\pi n)\right] =  {1\over 12} \left[  {3\ \Gamma^4 (1/4) \over 16\pi^3}  - 1 \right],\eqno(2.62)$$

$$\sum_{n=1}^\infty  n \left[1-  \tanh(\pi n/2 )\right] =   {1\over 12} \left[  {3  \Gamma^4 (1/4) \over 8 \pi^3}  - 1 \right],\eqno(2.63)$$

$$\sum_{n=1}^\infty  n \left[ \coth(\pi n) -1\right] =   {1\over 4} \left[ {1\over 3}  - {1 \over \pi}\right],\eqno(2.64)$$

$$\sum_{n=1}^\infty  n \left[ \coth(2\pi n) -1\right] =   {1\over 4}\left[ {1\over 3}  - {1 \over 2 \pi} \left[ 1+   { \Gamma^4 (1/4) \over 16 \pi^2}\right]\right],\eqno(2.65)$$

$$\sum_{n=1}^\infty   {n^2 \over  \sinh^2(\pi n)} =  { 1 \over 8  \pi^2} \left[ {\Gamma^8 (1/4) \over  192  \pi^4} - 1   \right],\eqno(2.66)$$

$$\sum_{n=1}^\infty   {n^2 \over  \sinh^2(2\pi n)} =  {1 \over 32 \pi^2}  \left[  {\Gamma^8 (1/4) \over  96  \pi^4} -  {\Gamma^4 (1/4) \over  8 \pi^2}  -1  \right],\eqno(2.67)$$

$$\sum_{n=1}^\infty  { n^2 \over  \cosh^2 (\pi n)} = {\Gamma^4 (1/4)   \over 64  \pi^4}  \left[ 1-    {\Gamma^4 (1/4)   \over 24  \pi^2}\right].\eqno(2.68)$$

\section{Series of $\hbox{cosech} (\pi nx)$ }

Finally, let us investigate  the following series

$$S(x)=  \sum_{n=1}^\infty  {1\over \sinh(\pi n x)},\quad x >0. \eqno(3.1)$$
It seems that (3.1) is much simpler that series considered in the previous section.    However, its calculation is quite a difficult task.  In this section we will deduce a second kind singular integral equation, involving the Hilbert transform [10] and whose solution is related to (3.1).  Indeed, let  represent it in a different form first.  To do this,  we appeal to   integral representation (1.16) for the Riemann zeta- function,  substituting  it in (3.1) and interchanging  the order of integration and summation via the absolute and uniform convergence by the same arguments as above.  Hence employing  the series for the Riemann zeta-function,  we find

$$S(x)= {1\over \pi i} \int_{\mu -i\infty}^{\mu +i\infty}  (2^s-1)  \Gamma(s)  \zeta^2(s)  (2\pi x)^{-s} ds,\  \mu  >1.\eqno(3.2)$$
Meanwhile, the square of the Riemann zeta-function is represented by another  Ramanujan's  identity [13]

$$\zeta^2(s)= \sum_{n=1}^\infty {d(n)\over n^s},\  {\rm Re}  s  > 1,\eqno(3.3)$$
where $d(n)$  is the Dirichlet divisor function, i.e. the number of divisors of $n$, including 1 and $n$  itself (see [1]).  Hence from (3.1), (3.2) and straightforward calculations we find 

$$S(x) =    2 \sum_{n=1}^\infty d(n)  \left[ e^{-\pi n x} -   e^{- 2 \pi n x} \right],\  x >0.\eqno(3.4)$$
In the meantime,  in [12]  we established the following particular case of the Voronoi summation formula, related to (3.4), namely

$$\sum_{n=1}^\infty d(n)  \left[ e^{-\pi n x}     + {2\over \pi x} \left( e^ {- 4\pi n/x} \hbox{Ei}\left({4\pi n\over x} \right) + 
   e^ {4\pi n/x} \hbox{Ei}\left(- \  {  4\pi n\over x} \right)\right)\right] = {1\over 4} + {\gamma - \log(\pi x) \over \pi x }, \eqno(3.5)$$
where $x >0,\ \gamma$ is Euler's constant and 

$$\hbox{Ei} (z)= \int_{-\infty}^z {e^t\over t} dt$$
is the integral exponential function.  Hence  with simple substitutions (3.5) can be written in the form

$$\sum_{n=1}^\infty d(n)  \left[ e^{-\pi n x}     + {4\over \pi x} \int_0^\infty { e^ {- 4\pi n t /x}\   t  dt \over 1- t^2} \right] = {1\over 4} + {\gamma - \log(\pi x) \over \pi x }, \eqno(3.6)$$ 
where,  as usual,  the integral in the neighborhood of $t=1$ is understood in the Cauchy principal values sense.  So, combining with (3.4), we find 

 $$ S(x) -  {4 \over \pi x} \sum_{n=1}^\infty d(n)  \int_0^\infty { \left( e^ {- 2\pi n t /x} - e^ {- 4\pi n t /x} \right) \   t  dt \over 1- t^2}  +  {4 \over \pi x} \sum_{n=1}^\infty d(n)  \int_0^\infty { e^ {- 4\pi n t /x} \   t  dt \over 1- t^2}$$

 $$=  {\gamma -  \log(\pi x/2) \over \pi x }.\eqno(3.7)$$
In order to change the order of integration and summation  in (3.7), we use the asymptotic behavior of the arithmetic function  $d(n)= O(n^\varepsilon)$ for all $\varepsilon > 0$  [1]  and split each  integral on four  integrals over the intervals, containing  as end-points $t=0,  \infty$ and $1\pm \delta,\  \delta >0.$   Let us show this on the integral

$$P.V. \int_0^\infty { e^ {- 4\pi n t /x}\   t  dt \over 1- t^2} = \int_0^{1/2} { e^ {- 4\pi n t /x}\   t  dt \over 1- t^2}  +
\lim_{\delta \to 0} \left( \int_{1/2}^{1-\delta} +  \int_{1+\delta}^2  \right) { e^ {- 4\pi n t /x}\   t  dt \over 1- t^2}$$

 $$+  \int_2^\infty { e^ {- 4\pi n t /x}\   t  dt \over 1- t^2}.$$
 Hence for any fixed $x >0$ we have the estimates

 $$ \int_0^{1/2} { e^ {- 4\pi n t /x}\   t  dt \over 1- t^2} \le \   8 x^3   \int_0^{1/2} {  t  dt \over 6 x^3  +  (4\pi n t)^3}$$

$$ \le { (x/\pi)^{3/2}  \over  2 \sqrt 6 \   n^{3/2} } \int_0^{1/2} {  dt \over \sqrt t} =  { (x/\pi)^{3/2}  \over  \sqrt {12} \   n^{3/2} }, $$

$$  \int_2^\infty { e^ {- 4\pi n t /x}\   t  dt \over t^2-1} \le   {1\over 8}   \left( {x\over \pi n}\right)^2  \int_2^\infty { dt \over  t (t^2-1)  } = {C_x\over n^2}. $$
 Concerning the middle  integrals, we appeal to the Lagrange theorem to write 

 $$\left| P.V. \int_{1/2}^2  { e^ {- 4\pi n t /x}\   t  dt \over 1- t^2} \right| =  \left| \int_{1/2}^2  {\left( e^ {- 4\pi n t /x} -   e^ {- 4\pi n /x}\right) \   t  dt \over 1- t^2} \right| +    e^ {- 4\pi n /x} \left| P.V. \int_{1/2}^2  {   t  dt \over 1- t^2}\right| $$

$$\le   {4\pi n\over x}\   e^ {- 2\pi n /x}  \int_{1/2}^2  {t  dt \over 1+ t} +  
  e^ {- 4\pi n /x} \left| \lim_{\delta \to 0}  \left(  \int_{1/2}^{1-\delta}  { t  dt \over 1- t^2} - \int_{1+\delta}^2  { t  dt \over  t^2-1}\right) \right|$$

 $$=    {e^ {- 2\pi n /x} \over x} \left[ 2\pi n (3- 2\log2)   +    e^ {- 2\pi n /x} \log 2\right].$$
Consequently, owing to these estimates the desired interchange is guaranteed,  taking $\varepsilon \in (0, 1/2)$.  Hence recalling (3.4),  equality (3.7) becomes 

 $$  S(x) -  {2 \over \pi x} \int_0^\infty  S\left( {2 t \over x}\right) \   {t  dt \over 1- t^2}  +  {4 \over \pi x}  \int_0^\infty  \sum_{n=1}^\infty d(n)  { e^ {- 4\pi n t /x} \   t  dt \over 1- t^2} =  {\gamma -  \log(\pi x/2) \over \pi x }.\eqno(3.8)$$ 
Letting in (3.8) $2x$ instead of $x$ and multiplying by $2$ the obtained equation, we subtract it from (3.8) with the use of (3.4), ending up with the following integro-functional equation 

  $$ S(x)- 2 S(2x)  +  {2 \over \pi x} \int_0^\infty  \left[   S\left( { t \over x}\right)-  2 S\left( {2 t \over x}\right) \right]  \   {t  dt \over 1- t^2}  =  { \log 2 \over  \pi x },$$ 
 which is up to simple substitutions coincides with the second kind singular integral equation,  associated with the Hilbert transform  of   $f(x)=  S(x)- 2 S(2x)$ 

  $$  f (x)  +  {2 \over \pi } \int_0^\infty    {  f( t)   x t \  dt \over 1- (xt)^2}  =  { \log 2 \over  \pi x},\quad x >0.\eqno(3.9)$$ 
This is in fact an exceptional  case of the Fox second kind integral equation (see [10], Section 11.15), whose solution cannot be written  as in [10] in terms of the inverse Mellin transform.   Nevertheless, the homogeneous equation (3.9) is solved recently by the author  (see details in [14], Corollary 2), using the method of compositions with the Fourier and Hartley transforms.   Now, defining the Mellin transform of $f$, for instance,   in  $L_2(\mathbb{R}_+)$ [10] as 

$$F(s)=  \int_0^\infty f(x) x^{s-1} dx,\ s \in \sigma = \left\{s \in \mathbb{C},\  {\rm Re} s = {1\over 2}\right\},$$
and taking into account the known formula 

$${1\over \pi} PV \int_0^\infty {t^{s-1} \over 1-t } =  \cot(\pi s),\   0 < {\rm Re s} < 1,\eqno(3.10)$$
the Mellin transform of the left-hand side of (3.9) is equal to

$$F(s) -   F(1-s) \tan\left({\pi s\over 2}\right).$$
But

$$F(s)=  S^*(s) \left(1 - 2^{1-s} \right),$$
where $S^*(s)$ is the Mellin transform of $S(x)$.  Recalling (3.2) and the functional equation (1.21) for the Riemann zeta-function one easily verifies that 

 $$F(s) -   F(1-s) \tan\left({\pi s\over 2}\right)=0,$$
i.e.  $f$ satisfies the homogeneous equation (3.9), which has an infinite number of solutions as it is shown in [14].   We can also get it directly, applying the Hilbert transform to both sides of (3.9) and taking into account that integral (3.10) is zero for $s=1/2$. 

Finally, we observe from (1.19),  (1.20), (3.2) that the same scheme can be applied to the series (1.2) with $\alpha=0$,
because they relate to (3.1), for instance,  via identities 

$$\sum_{n=1}^\infty  {1\over \sinh(\pi n x)} =   \sum_{n=1}^\infty   \left[\coth(\pi n x/2) -1 \right]  -  \sum_{n=1}^\infty   \left[\coth(\pi n x) -1 \right],\ x >0,$$

$$  \sum_{n=1}^\infty   \left[\tanh(\pi n x/2) -1 \right] =  \sum_{n=1}^\infty   \left[\coth(\pi n x) -1 \right] - \sum_{n=1}^\infty  {1\over \sinh(\pi n x)},\ x >0.$$

\bigskip
\centerline{{\bf Acknowledgments}}
\bigskip

The work was partially supported by CMUP (UID/MAT/00144/2013), which is funded by FCT(Portugal) with national (MEC) and European structural funds through the programs FEDER, under the partnership agreement PT2020. The author is sincerely indebted to the referee for pointing out the pioneer references on the topic and his comment, which led to the idea to extend the method for  series involving the hyperbolic tangent and cotangent functions.  

\bigskip
\centerline{{\bf References}}
\bigskip
\baselineskip=12pt
\medskip
\begin{enumerate}

\item[{\bf 1.}\ ]    T. M. Apostol, {\it Modular Functions and Dirichlet Series in  Number Theory},  2nd ed. Springer,  New York (1990).

\item[{\bf 2.}\ ]  J.M. Borwein and P.B. Borwein,  {\it   Pi and  the AGM: A Study in Analytic Number Theory and Computational Complexity},  Wiley,  New York (1987).

\item[{\bf 3.}\ ]  J.M. Borwein and I.J. Zucker,  Elliptic integral evaluation of the Gamma-function at rational values of small denominators, {\it  IMA J. Numerical Analysis}, \ {\bf 12} (1992),  519-526.

\item[{\bf 4.}\ ]  M.L. Glasser, V.G. Papageorgiou and T.C. Bountis,  Melnikov's function for two-dimensional mappings, {\it  SIAM J. Appl. Math.}, \ {\bf 49} (1989),  N 3,  692-703.

\item[{\bf 5.}\ ]  Ya. Komori,  K. Matsumoto  and H. Tsumura,  Infinite series involving hyperbolic functions, {\it  Lithuanian Math. J.}, \ {\bf 55}, N 1  (2015),  102- 118.

\item[{\bf 6.}\ ]   C.B. Ling, On summation of series of hyperbolic functions.  {\it SIAM J. Math. Anal.}  {\bf 5},   (1974),  551-562.

\item[{\bf 7.}\ ]   A. Erd\'elyi, W. Magnus, F. Oberhettinger and F.G. Tricomi, {\it Higher Transcendental Functions}, Vols  I, II and III, McGraw-Hill, New York, London and Toronto (1953).

\item[{\bf 8.}\ ]  C. Nasim, A summation formula involving $\sigma(n)$, {\it Transactions of the American Mathematical Society}, \ {\bf 192} (1974), 307- 317.

\item[{\bf 9.}\ ]A.P. Prudnikov, Yu. A. Brychkov and O. I. Marichev, {\it Integrals and Series:} Vol. I: {\it Elementary
Functions}, Gordon and Breach, New York (1986);    Vol. II:  {\it Special functions},  Gordon and Breach, New York (1986);   Vol. III:  {\it More special functions},   Gordon and Breach, New York (1990).

\item[{\bf 10.}\ ]  E.C. Titchmarsh, {\it  An Introduction to the Theory of Fourier Integrals},    Chelsea, New York  ( 1986).

\item[{\bf 11.}\ ]  E.T. Whittaker and  G.N. Watson  {\it A Course in Modern Analysis}, 4th ed. Cambridge University Press, Cambridge (1990) .

\item[{\bf 12.}\ ]  S. Yakubovich,  A general class of Voronoi's and Koshlyakov- Ramanujan's summation formulas involving $d_k(n)$,  {\it Integral Transforms and Special Functions},  {\bf 22}, N 11  (2011),  801-821. 

\item[{\bf 13.}\ ]  S. Yakubovich,  Integral and series transformations via Ramanujan's identities and Salem's type equivalences to the Riemann hypothesis,  {\it Integral Transforms and Special Functions},  {\bf 25}, N 4  (2014),  255-271.

\item[{\bf 14.}\ ]  S. Yakubovich, On the half-Hartley transform, its iteration and compositions with Fourier transforms. {\it J. Integral Equations Appl.}  {\bf 26}, N 4  (2014),  581-608.

\item[{\bf 15.}\ ]  I.J. Zucker,  The summation of series of hyperbolic functions.  {\it SIAM J. Math. Anal.}  {\bf 10}, N 1  (1979),  192- 206.

\end{enumerate}

\vspace{5mm}

\noindent S.Yakubovich\\
Department of Mathematics,\\
Faculty of Sciences,\\
University of Porto,\\
Campo Alegre st., 687\\
4169-007 Porto\\
Portugal\\
E-Mail: syakubov@fc.up.pt\\

\end{document}